\documentclass[12pt,leqno]{article}

\newlength{\noteWidth}
\long\def\notes#1{\ifinner
             {\tiny #1}
             \else
              \marginpar{\parbox[t]{\noteWidth}{\raggedright\tiny #1}}
               \fi}
\def\notes#1{}

\usepackage{latexsym}
\usepackage{graphics}
\usepackage{amsmath}
\usepackage{xspace}
\usepackage{amssymb}
\usepackage{latexsym,fullpage}
\usepackage{amsmath,amsthm}
\usepackage{epsfig}
\usepackage{pst-all}

\usepackage[usenames]{color}

\definecolor{MyDarkBlue}{rgb}{0,0.008,0.25}
\definecolor{MyBlueGreen}{rgb}{.015,0.4,0.5}
\definecolor{MyDarkRed}{rgb}{.5,0.05,0.05}

\def\ignore#1{{}}

\newcommand {\be}[1]{\begin{equation}\label{#1}}
\newcommand {\ee}{\end{equation}}
\newcommand {\bea}{\begin{eqnarray}}
\newcommand {\eea}{\end{eqnarray}}

\newcommand{\Exp}{\ensuremath{\operatorname{Exp}}\xspace}

\newcommand{\pr}{\mathbb{P}}
\newcommand{\E}{\mathbb{E}}

\newtheorem{theorem}{Theorem}[section]
\newtheorem{prop}[theorem]{Proposition}
\newtheorem{lemma}[theorem]{Lemma}

\def\barT{{\bar{T}}}

\def\varble{\,\cdot\,}

 \def\transpose{{\hbox{\rm\tiny T}}}

\def\clL{{\cal L}}

\def\clS{{\cal S}}

\def\breakMed{\\[.2cm]}

\newcommand{\field}[1]{\mathbb{#1}}

\def\Re{\field{R}}

\def\nat{\field{Z}_+}

\def\bfmath#1{{\mathchoice{\mbox{\boldmath$#1$}}%
{\mbox{\boldmath$#1$}}%
{\mbox{\boldmath$\scriptstyle#1$}}%
{\mbox{\boldmath$\scriptscriptstyle#1$}}}}

\def\bfmq{\bfmath{q}}

\def\bfmz{\bfmath{z}}

\def\bfmA{\bfmath{A}}

\def\bfmD{\bfmath{D}}
\def\bfmL{\bfmath{L}}

\def\bfmQ{\bfmath{Q}}

\def\bfmS{\bfmath{S}}

\def\bfmX{\bfmath{X}}

\def\bfmZ{\bfmath{Z}}

\def\bfmx{\bfmath{x}}

\newcounter{rmnum}
\newenvironment{romannum}{\begin{list}{{\upshape (\roman{rmnum})}}{\usecounter{rmnum}
\setlength{\leftmargin}{24pt}
\setlength{\rightmargin}{16pt}
\setlength{\itemindent}{-1pt}
}}{\end{list}}

\newcounter{anum}

\def\Lemma#1{Lemma~\ref{#1}}
\def\Proposition#1{Proposition~\ref{#1}}
\def\Theorem#1{Theorem~\ref{#1}}
\def\Section#1{Section~\ref{#1}}
\def\Figure#1{Figure~\ref{#1}}

\def\FRAC#1#2#3{\genfrac{}{}{}{#1}{#2}{#3}}

\def\ddt{{\mathchoice{\FRAC{1}{d}{dt}}%
{\FRAC{1}{d}{dt}}%
{\FRAC{3}{d}{dt}}%
{\FRAC{3}{d}{dt}}}}

\def\eqdef{\mathbin{:=}}

\def\clC{{\cal C}}

\def\clF{{\cal F}}

\title{On Exponential Ergodicity \\
of Multiclass Queueing Networks}

\author{{\sf David Gamarnik}\thanks{MIT Sloan School of Management, Cambridge, MA,  02139 e-mail: {\tt
gamarnik@mit.edu}}\\ MIT
\and
 {\sf Sean Meyn}\thanks{Department of Electrical and Computer Engineering and the Coordinated Sciences Laboratory, University
of Illinois at Urbana-Champaign, Urbana, IL 61801, e-mail: {\tt meyn@control.csl.uiuc.edu}}\\
UIUC}

\begin{document}

\maketitle


\begin{abstract}
One of the key performance measures in queueing systems is the
exponential decay rate of the steady-state tail probabilities of the
queue lengths.  It is known that if a corresponding fluid model is
stable and the stochastic primitives have finite moments, then  the
queue lengths also have finite moments, so that the tail probability
$\pr(\cdot >s)$ decays faster than $s^{-n}$ for any $n$.   It is
natural to conjecture that the decay rate is in fact exponential.

In this paper an example is constructed to demonstrate that   this
conjecture is false.  For a specific stationary policy applied to a
network with exponentially distributed interarrival and service
times it is shown that the corresponding fluid limit model is stable,
but the tail probability for the buffer length   decays slower than
$s^{-\log s}$.
\end{abstract}



\section{Introduction.}
A key performance measures in queueing models is   the decay rate of
the queue length distribution in steady state
\cite{glywhi94b,ganocowis04,dufoco95,CTCN}.
Except in trivial cases,  the
decay rate is at best exponential for networks with a finite number
of servers. Moreover,  an exponential decay rate can be verified in
many queueing models either by direct probabilistic arguments such
as Kingman's classical bound for the G/G/1 queueing system, using
Lyapunov function type arguments
\cite{mey05b,kummey96a,GamarnikZeevi} (see in particular
Sections~16.3 and 16.4 of \cite{metwee}),   using large deviations
techniques \cite{large_deviations}, and even specialized techniques
based on a fluid limit model \cite[Theorem~4]{mey01a}.

\begin{figure}\label{figure:RybkoStolyar}
\begin{center}
\includegraphics[width=.75\hsize]{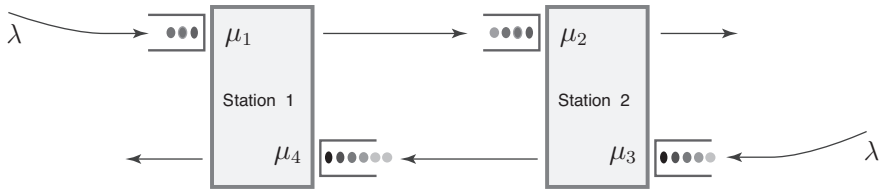}
\end{center}
\caption{Rybko-Stolyar network}
\end{figure}

Except in special cases, such  as single-class queueing networks
\cite{GamarnikZeevi}, verifying the existence of an exponential tail
is non-trivial precisely because verifying stability of a
multi-class network is no longer straightforward.     It is now well
known that the standard load condition $\rho_\bullet<1$ is not
sufficient for stability.     This was demonstrated in the two
seminal papers of Kumar and Seidman~\cite{kumsei90a} and Rybko and Stolyar \cite{rs}, based on the network depicted in \Figure{figure:RybkoStolyar}, henceforth called the \textit{KSRS network}.

Motivated in part by these examples, over the past ten years new
methods have been developed to verify stability.  The most general
techniques are based on the fluid limit model (recalled below
equation \eqref{scaleQ}) starting with the work of  Maly{\v{s}}ev
and Men{$'$}{\v{s}}ikov~\cite{malmen79} and Stolyar~\cite{stolyar}.
This was extended to a broad class of multiclass networks by Dai
\cite{dai}.

Dai showed that stability of the fluid limit model together with some mild conditions on the network implies positive Harris recurrence for a Markovian state process, which implies in particular the existence of a unique steady state for the underlying queueing network.  This result was extended in Dai and Meyn~\cite{DaiMeyn96} where it is shown that the queue lengths have finite moments in steady state up to order $p$ if the   stochastic primitives of the network (interarrival and service
time distributions) have finite moments up to $p+1$.

As a direct implication of the main result of \cite{DaiMeyn96},  if the stochastic primitives have an exponential moment, and if the fluid limit  model is stable, then
the decay rate of the queue length  distribution in steady state is faster than than any polynomial, in the sense that, for any $p\ge 1$,
\begin{equation}
\lim_{s\to\infty}  r(s)  \pr\{ \|Q(0)\| \ge s\}  =0, \label{tail}
\end{equation}
where $r(s)= s^p$, and the probability is with respect to a stationary version of the queue.

It is then natural to conjecture that an exponential bound holds, so
that \eqref{tail} holds with $r(s)=e^{\theta s}$ for  some
$\theta>0$, provided the stochastic primitives possess
exponentially decaying tails.    The purpose of this paper is to
construct a particular  stationary policy and a particular network
to demonstrate that this conclusion in fact does not hold. \notes{spm: the conjecture can be made, but the conclusion is false;)}
Moreover,  a polynomial rate is about the best that can be attained:
In the example it is found that \eqref{tail} cannot hold for any
sequence $\{r(s)\}$ satisfying $\liminf_{s\to\infty}
r(s)/s^{\log(s)}>0$.

The example is a particular instance of the KSRS model in which the interarrival and service times have exponential distributions.   The scheduling policy is stationary, and the corresponding fluid limit model is stable, yet the queue length process has heavy tails in steady-state.
The policy is based on carefully randomizing
between a stable policy, and the unstable policy introduced by Rybko and Stolyar.

\medskip

The remainder of the paper is organized as follows. In the following section we describe the model and the main result. In Section~\ref{stable} we present a short proof of  the earlier result  of \cite{DaiMeyn96} in the special setting of this paper.     The construction of the counterexample is provided in Section~\ref{section:network}.
The details of the proof of the main result are contained
in Sections~\ref{section:stability} and \ref{section:instability}.
Some technical arguments are placed in an appendix.

We close the  introduction with some notational conventions: The
$i$th unit vector in $\Re^N$ is denoted $e_i$ for $1\le i\le N$. The
$L_1$-norm, always denoted $\|\cdot\|$, is defined by
$\|a\|=\sum_{1\leq i\leq N}|a_i|$ for $a\in\Re^N$. Given two
positive real valued functions $f(x),g(x)$, the notation $f=O(g)$
stands for $f(x)\leq C g(x)$ for some constant $C>0$ and all $x\geq
0$.  Similarly, $f=\Omega(g)$ means $f(x)\geq Cg(x)$ for all $x\geq
0$ and $f=\Theta(g)$ means $f=O(g)$ and $f=\Omega(g)$ at the same
time.
Given a Markov chain or a Markov process $\bfmZ$
defined on a space ${\cal Z}$ and given a probability measure $\nu$
on ${\cal Z}$, we let $\pr_\nu(Z(t))$ denote the law of $Z(t)$
initialized by have $Z(0)$ distributed according to $\nu$.
Specifically, for every $x\in {\cal Z}$, $\pr_x(Z(t))$ is the law of
$Z(t)$ conditioned on $Z(0)=x$. The notations
$\E_\nu[\cdot],\E_x[\cdot]$ have similar meaning corresponding to
the expectation operator.

\section{Model description and main result}
\label{main}

The model and the definitions of this paper follow closely those
of \cite{DaiMeyn96}. We consider a multiclass queueing network
consisting of $J$ servers denoted simply by $1,2,\ldots,J$. Each
customer class is associated with an exogenous arrival process which
is assumed to be a renewal process with  rate $\lambda_i$.
Here $N$ denotes the number of classes.
It is possible that $\lambda_i=0$ for some of the classes, namely, no
external arrival is associated with this class. We let
$\lambda=(\lambda_i), 1\leq i\leq N$. Each server is unit speed and
can serve customers from a fixed set of customer classes
$i=1,2,\ldots,N$. The classes are associated with servers using the
constituency matrix $C$, where $C_{ij}=1$ if class $i$ is served by
server $j$, and $C_{ij}=0$ otherwise. It is assumed that each class
is served by exactly one server, but the same server can be
associated with many classes. Each class is associated with a
buffer at which the jobs are queued waiting for service.  The queue
length corresponding to the jobs in buffer $i$ at time $t$ is
denoted by $Q_i(t)$, and $Q(t)$ denotes the corresponding $N$-dimensional buffer.

It is assumed that routing is deterministic:
The routing matrix $R$ has entries equal to zero or one.
Upon service completion, a job in class $i$ proceeds to buffer $i_+$, where $i_+$ denotes the
index satisfying $R_{i\, i_+}=1$, provided such an index exist.
In this case buffer~$i$ is called an \textit{internal buffer}.
If no such index exists, then this is an \textit{exit buffer},
and the completed job leaves the network.  The network is assumed to be open, so that $R^N= 0$.

For the purpose of building a counterexample we
restrict to a Markovian model:  The arrival processes are assumed Poisson, and the
service times have exponential distributions.  It is also convenient to \textit{relax} the assumptions above and allow infinite rates for service at certain queues.  That is, the corresponding service time is \textit{zero}.   For this reason it is necessary to show that the main result of \cite{DaiMeyn96} can be extended to this setting.

We can express the evolution of the queue length process as
\begin{equation}
\begin{aligned}
Q(t)  &= Q(0)
+ \sum_{i,j=1}^N   [-e_i+R_{ij}] D_i(t)  +  A(t)
\\
 &= Q(0) +    [-I+R^\transpose] D(t)  +  A(t),
\end{aligned}
\label{RoutingModelD}
\end{equation}
where $\bfmA_i$ is the cumulative arrival process to buffer $i$, and $\bfmD_i$ is the cumulative departure process from buffer $i$.
The second equation is in vector form with  $R$ equal to the routing matrix, and $D(t),A(t)$ the $N$-dimensional vectors of departures and arrivals.
\notes{DG:  I am confused about the first equation. Is it vector form?
\\
spm:  ok now?}

If the service rate $\mu_i$ is finite, then the departure process can be expressed,
\[
 D_i(t) = S_i(Z_i(t)) ,\qquad i=1,\dots, N,\  t\ge 0,
\]
where $\bfmS_i$ is a Poisson process with rate $\mu_i$, and $Z_i(t)$
is the cumulative busy time at buffer~$i$.  All of these processes
are assumed right continuous.

The following  assumptions are imposed on the policy that determines
$\bfmZ$:   It is assumed throughout the paper that the policy is
\textit{stationary}.   In this setting, for buffers with finite
service rate this means that $\ddt Z_i(t)$ is piecewise constant,
and when the derivative exists it can be expressed as a fixed
function of the queue-length process, $\ddt Z_i(t)=\phi_i(Q(t))$.
Moreover, for each $i$ satisfying $\mu_i=\infty$ there is a fixed
set of values $\Xi_i \subset\nat^N$ such that the contents of buffer
$i$ are drained at the instant $Q(t)\in \Xi_i$.

We  also considered randomized stationary policies.  In this case
$\phi_i$ is a randomized function of $Q(t)$, and the draining of a
buffer with infinite service rate occurs with some probability
depending upon the particular value of $Q(t)\in \Xi_i$ observed at
time $t$.

It is assumed that the policy is \textit{non-idling}:
For each station we have $\sum \ddt Z_i(t) =1$ when $\sum Q_i(t)>0$,
where the sum is over $i$ at the given station.


Throughout much of the paper we restrict to the KSRS model in which the routing and constituency matrix are expressed,
\[
\textstyle
R^\transpose=\begin{bmatrix} 0 & 0 & 0 &0\\
                                        1 & 0 & 0 & 0\\
                                          0 & 0 & 0&0\\
                                          0 & 0 & 1&0
\end{bmatrix}
\qquad
C=\begin{bmatrix}    1 &0&0& 1\\
                                          0 & 1 &1& 0
\end{bmatrix}
\]
The network is symmetric:  The two non-null arrival processes are independent Poisson processes with rates $\lambda_2=\lambda_4>0$, and the service rates at buffers~2 and 4 are finite and equal, $\mu_2=\mu_4<\infty$.  The service rates at buffers~1 and 3 are infinite.  Hence, for example, if at any moment priority is given to buffer~1, then all of the contents pass instantly to buffer~2.


The transition semigroup for $\bfmQ$ is  denoted,
\[
P^t(x,A) = \pr \{Q(t)\in A \mid Q(0) =x\},\qquad x\in\nat^n,\ t\ge
0, A\subset\nat^N,
\]
and $\pr_x(\varble)$ the probability law corresponding to the
initial state $Q(0)=x$.   A probability measure $\pi$ on the state
space $\nat^N$ is invariant if $\bfmQ$ is a stationary process when
initialized using $\pi$.
This is equivalently expressed by the invariance
equations,
\[
\pi(y) = \sum_x\pi(x)P^t(x,y) ,\qquad y\in\nat^N,\ t\ge 0 .
\]

We say that $\pi$ has an \textit{exponential tail} if for some
$\theta>0$,
\[
\sum_x \pi(x)e^{\theta \|x\|} < \infty.
\]
The Markov process is called \textit{exponentially ergodic} if $\pi$ exists, and for some $\theta>0$ and each   $x,y$,
\begin{equation}
\lim_{t\to\infty} e^{\theta t} | \pr_x(Q(t)=y) -\pi(y)|  =0.
\label{geo}
\end{equation}
Exponential ergodicity implies an exponential tail.  The proof of \Proposition{geotail} is provided in \Section{stable}.
\begin{prop}
\label{geotail}
Consider the network \eqref{RoutingModelD} in which $\mu_i<\infty$ for each exit buffer.  Assume that the network is controlled using a stationary policy.  If $\bfmQ$ is exponentially ergodic, then it has an exponential tail.
\end{prop}

We now   construct the \emph{fluid limit model} associated with the
network \eqref{RoutingModelD}.  To
emphasize the dependence on the initial state $Q(0)=x$ we denote the queue and allocation trajectories by $Q(\cdot,x),Z(\cdot,x)$.  The scaled initial condition is defined for $\kappa>0$ by,
 \begin{equation}
 x^\kappa \eqdef \frac{1}{\kappa}\lfloor \kappa x\rfloor,\qquad x\in\Re_+^N,
\label{xr}
\end{equation}
and the scaled processes are defined for $t\ge 0$ via,
\begin{equation}
q^\kappa(t;x^\kappa)  \eqdef  \frac{1}{\kappa} Q(\kappa t; \kappa
x^\kappa), \quad z^\kappa(t;x^\kappa)  \eqdef \frac{1}{\kappa}
Z(\kappa t; \kappa x^\kappa). \label{scaleQ}
\end{equation}
Observe that $x^\kappa\in\Re_+^N$ satisfies $\kappa
x^\kappa\in\nat^N$ for each $\kappa$. For each $x\in\nat^N$ and
$\omega\in\Omega$ we let $\clL_x(\omega)$ denote the set of all
possible fluid limits,
\[
\clL_x(\omega) = \Bigl\{\  \parbox{.63\hsize}{u.o.c.\ subsequential limits of  $\{q^\kappa(t;x^\kappa,\omega) ,\ z^\kappa(t;x^\kappa,\omega)\}$\ } \Bigr\},
\]
where `u.o.c.' means that  convergence is uniform on compact time
intervals as $\kappa_i\to\infty$ for some subsequence
$\{\kappa_i\}$. The \textit{fluid limit model} is   the union
$\displaystyle\clL \eqdef \cup_{x\in\nat^n} \clL_x$.  It is clear
that any fluid limit must satisfy the fluid model equation,
\begin{equation}
q(t)  = q(0) +    [-I+R^\transpose] z(t)  +  \lambda t,
\label{RoutingModelDfluid}
\end{equation}
in which $(\bfmz,\bfmq)$ satisfy assumptions analogous to $(\bfmZ,\bfmQ)$ \cite{dai,CTCN}.
\notes{DG: shall we say "fluid equation of the form ..." plus citation.
Otherwise, what is $q,B,z$?
\\
spm: o.k.?  Let me know if I'm overdoing the ctcn references!}

The fluid limit model $\clL$ is said to be \textit{stable}  if there
exists $\Omega_0\subset \Omega$ satisfying $\pr\{\Omega_0\} =1$, and
$T_0>0$ such that $q(t) = 0$ whenever $t\ge T_0$,
$\omega\in\Omega_0$,  $\bfmq\in\clL(\omega)$ and $\|q(0)\| =1$.

Naturally, the fluid limit model as well as the conditions for
stability may depend on the scheduling policy.  In many cases (such
as the G/G/1 queue)  the set of fluid limits $\clL_x$ is a
deterministic singleton for each $x$.

The following is the key motivating result for our paper.
\begin{theorem}
\label{theorem:DaiMeyn}
Consider the network model \eqref{RoutingModelD}
controlled using a stationary non-idling policy.   Suppose that $\mu_i<\infty$ at each exit buffer, and that the fluid limit
model is stable. Then,
\begin{romannum}
\item
 {\em
$\bfmQ$ is aperiodic and positive Harris recurrent}:  There is a
unique invariant measure $\pi$ such that the distributions converge in total variation norm for each initial condition $x\in\nat^N$,
\[
\lim_{t\to\infty} \bigl( \sup_{y\in\nat^N}| \pr_x(Q(t)=y)
-\pi(y)|\bigr)  =0.
\]
\item
 {\em
The invariant measure has polynomial moments}:
For each $p\ge 1$,
\begin{align}
\sum_{x\in\nat^N}\pi(x) \|x\|^p<\infty . \label{eq:momentsfinite}
\end{align}
\end{romannum}
\end{theorem}

Theorem~\ref{theorem:DaiMeyn} asserts that the model \eqref{RoutingModelD} is positive Harris recurrent with polynomial moments of order $p$ for \textit{every integer $p$} when the interarrival and service times are exponentially distributed.
This suggests that $\pi$ will have an exponential tail, so that $\E_{\pi}[e^{\theta'\|Q(0)\|}]<\infty$ for some  $\theta'>0$. We now show that this conjecture does not hold
true.

Theorem 4.1  of \cite{DaiMeyn96} considers general  models with
renewal inputs, and also establishes rates of convergence to
stationarity and other ergodic theorems for the model.   In
Theorem~\ref{theorem:DaiMeyn}  we have extracted the part that is most relevant to the counterexample described next.

\begin{theorem}
\label{theorem:MainResult}
Consider the Rybko-Stolyar model described by the Markovian model \eqref{RoutingModelD} as follows,
\begin{equation}
\begin{aligned}
Q_i(t)&=Q_i(0)+A_i(t)-D_i(t), ~~i=1,3 \\
Q_i(t)&=Q_i(0)+D_{i-1}(t)-D_i(t),~~i=2,4.
\end{aligned}
\label{RSQ}
\end{equation}
There exist network parameters and a stationary policy satisfying,
\begin{enumerate}
\item The interarrival and service times are mutually independent with exponential distribution.

\item The fluid limit model is stable (hence  there exists a unique invariant measure $\pi$.)

\item
The invariant measure $\pi$ satisfies
\begin{align}
\label{eq:expmomentInfinite}
\E_{\pi}[\Psi(\|Q(0)\|)]=\infty,
\end{align}
where  $\Psi(s)=s^{\log s}$ for $s>0$, and $\Psi(0)=0$.
In particular, the invariant measure $\pi$ does not have an exponential tail and the Markov process is not exponentially ergodic.
\end{enumerate}
\end{theorem}

Theorem~\ref{theorem:MainResult}  establishes that the result
(\ref{eq:momentsfinite}) of Theorem~\ref{theorem:DaiMeyn} is nearly tight, modulo the $\log$ term in the exponent.


The proof of \Theorem{theorem:MainResult} is technical. It is
simplified substantially through our adoption of a  relaxation in
which the service rates at buffers 1 and 3 are infinite.     The
resulting process violates the assumptions of Dai and
Meyn~\cite{DaiMeyn96}, but we show in the following section that
these results carry over to this more general setting to yield
Theorem~\ref{theorem:DaiMeyn}.

\section{Stability of Markovian networks}
\label{stable}

In this section we establish some general properties of the  model
\eqref{RoutingModelD}.   In this section only  we consider the
embedded chain obtained via uniformization.  It is known that
geometric ergodicity\footnote{The definition of g.~ergodicity is precisely   \eqref{geo} in discrete time.}
of this chain is equivalent to  exponential
ergodicity of the process \cite{dowmeytwe95a}.

Uniformization must be applied with care when service rates can be
infinite.  We  define sampling times   $\{\tau_n, n\geq 0\}$
corresponding to jumps of a Poisson process derived from the
arrival-service process $\{\bfmA ,\bfmS\}$.  This is commonly
interpreted as sampling at  arrival epochs  and (real or virtual)
service completions. However, in sampling we restrict to non-zero
length service completions to avoid sampling twice at the same
instant!

Right-continuity implies that  $Q(\tau_n)=Q(\tau_n^+)$ for each $n$.

It is in fact simplest to ignore most of the network structure, and
consider a general Markov chain denoted $\bfmX$ on $\nat^N$
satisfying a version of the so-called skip-free property
\cite{metwee} along with a standard irreducibility condition.

The skip-free condition is defined with respect to the $L_1$ norm.
Under the assumption that the service rate at any exit buffer is finite
it follows that the increments $\big|\|X(t+1)\|-\|X(t)\|\big|$ are bounded
in the network model (\ref{RoutingModelD}) when $\{X(t)\eqdef Q(\tau_t) : t\in\nat \}$.  Note that the increments of the norm  $\|X(t+1) -X(t)\|$ are not bounded in a queueing model that allows
instantaneous transfer of buffer contents.

The following result establishes \Proposition{geotail}.
\begin{theorem}
\label{geotailRelax}
Suppose that $\bfmX$ is a Markov chain on $\nat^N$ satisfying the following
\begin{romannum}
\item
It is skip-free in the sense that
 for some constant $b_0$  and every initial condition,
\begin{equation}
-b_0\le  \|X(t+1)\| - \|X(t) \|\le b_0  ,\qquad t\ge 0. \label{skip}
\end{equation}

\item  The chain is  $0$-irreducible, in the sense that
\begin{equation}
\hbox{$\displaystyle \sum_t P^t(x,0)>0$ for each $x$.  } \label{irr}
\end{equation}

\item
$\bfmX$ is geometrically ergodic
\end{romannum}
 Then $\pi$ has an exponential tail.
\end{theorem}

\begin{proof}
The proof proceeds in two steps.

Denote the moment generating function by $h_\theta(x) = \E_x[e^{\theta \tau_0}]$, $\theta>0$,
$x\in\nat^n$, where $\tau_0$ is the first hitting time to the origin.
Theorem~15.2.4 and Theorem~16.3.2 of  \cite{metwee} imply that  $\pi(h_\theta)<\infty$ for sufficiently small $\theta>0$.

The second step is to compare $h_\theta$ with the function
$g_\theta(x) = e^{\theta\|x\|} $, $\theta>0$, $x\in\nat^N$. The
skip-free assumption implies the bound,
\begin{equation}
\bigl| \|X(t)\|- \|X(0)\| \bigr| \le b_0 t,\qquad t\ge 0,
\label{b0X}
\end{equation}
so that $\tau_0\ge b_0^{-1} \|x\|$ with probability one
when $X(0)=x$. Hence $g_{\theta/b_0}(x)\le h_\theta(x)$ for all $x$.
We conclude that,
\[
\E_\pi\bigl[ e^{ (\theta/b_0) \|X(0)\|} \bigr]
= \sum \pi(x) g_{\theta/b_0}(x)
\le \sum \pi(x) h_\theta(x)  <\infty.
\]
\end{proof}

The definition of the fluid limit model is defined for any Markov chain exactly as in the network model via \eqref{scaleQ}.
\begin{theorem}
\label{theorem:DaiMeynRelax}
Suppose that $\bfmX$ is a Markov chain on $\nat^N$ satisfying (i) and (ii) of
\Theorem{geotailRelax}.  Suppose moreover that   the fluid limit model is stable, in the sense that for some $T_0>0$ and a set
$\Omega_0\subset \Omega$ satisfying  $\pr\{\Omega_0\} =1$,
\notes{spm: trying to be super-precise here}
\[
\lim_{\kappa\to\infty}
  \frac{1}{\kappa} \| X(\kappa t; \kappa x^\kappa) \|= 0,\qquad \|x\|\le 1,\ t\ge T_0,\ \omega\in\Omega_0.
\]
Then $\bfmX$ is positive Harris  recurrent, and its unique invariant measure $\pi$ satisfies $\sum\pi(x) \|x\|^p <\infty $ for each $p\ge 1$.
\end{theorem}

\begin{proof}
In \Proposition{Drift} we establish a Lyapunov drift condition of the form:  For each $p\ge 1$ we can find a function $V$ and positive constants  $b$ and $\epsilon$ such that,
\begin{equation}
PV\, (x)  = \E[V(X(t+1)) \mid X(t)=x] \le V(x) - \epsilon
V^{1-\delta}(x)  + b\, ,
\label{MoulinesV3}
\end{equation}
where $\delta=(1+p)^{-1}$.  Moreoever, the function $V$ is equivalent to $\|x\|^{p+1}$:
\begin{equation}
 0<
  \liminf_{r\to\infty}  \Bigl( \inf_{  \|x\|=r} \frac{V(x)}{r^{p+1}}   \Bigr)
 \le
  \limsup_{r\to\infty}  \Bigl( \sup_{ \|x\|=r} \frac{V(x)}{r^{p+1}}   \Bigr)
  <\infty .
\label{polyGrowth}
\end{equation}
It follows from the Comparison Theorem of \cite{metwee}
that the steady-state mean of $V^{1-\delta}$ is finite, with the explicit bound $\pi(V^{1-\delta})\le b/\epsilon$.  This implies that the $p$th moment of $\bfmX$ is finite since $1-\delta=p/(p+1)$.
\end{proof}

The drift criterion \eqref{MoulinesV3} was introduced in the analysis of general state-space Markov chains in  \cite{douformousou04}.   Under this bound polynomial rates of convergence are obtained on the rate of convergence to steady state.  Using different methods, polynomial bounds on the steady-state buffer lengths and polynomial rates of convergence were obtained in \cite{DaiMeyn96} for stochastic networks based on a general result of \cite{tuotwe94a}.   The proof is simplified considerably in this countable state space setting.

To establish \eqref{MoulinesV3}   we take $V=V_p$ with,
\begin{equation}
V_p(x)= \E \Bigl[\sum_{t=0}^{\|x\|T}  \| X(t)\|^p   \Bigr],  \qquad \
x\in\nat^N, \label{V3fluidCRW}
\end{equation}
where   $T\ge 1$ is a sufficiently large fixed integer.
The growth bounds \eqref{polyGrowth} are suggested by the approximation,
\begin{equation}
\frac{1}{\kappa^{p+1}}V_p(\kappa x^\kappa) \approx
 \E \Bigl[\int_0^{\|x\|T}  \| x^\kappa(t;x^\kappa)\|^p \,  dt  \Bigr],\qquad x\in\Re_+^N,\ \kappa>0,
\label{polyGrowthB}
\end{equation}
with $\bfmx^\kappa$ defined as in \eqref{scaleQ} via $
x^\kappa(t;x^\kappa)  \eqdef  \kappa^{-1} X(\kappa t; \kappa
x^\kappa)$. \notes{Yikes!  We went from X back to Q suddenly.  All fixed now.
-spm}

\begin{prop}
\label{Drift}
The following hold under the assumptions of
\Theorem{theorem:DaiMeynRelax}:  for each $p=1,2,\dots$  the function $V$ defined in \eqref{V3fluidCRW} satisfies the drift condition
\eqref{MoulinesV3} and the bounds \eqref{polyGrowth}
with $\delta=(1+p)^{-1}$.
\end{prop}

To prove the proposition we first note that stability of the fluid model implies convergence in an $L_p$ sense.   The proof
of the almost sure limit in \Proposition{FLMstableB} uses equicontinuity of $\{ x^\kappa(t;x^\kappa) : t\ge 0 \}$, and the $L_p$ limit is obtained using the Dominated Converence Theorem since $\bfmL$ is a bounded sequence.  Recall that  $T_0$ and $\Omega_0$ are introduced in the definition of stability for the fluid limit model.
\begin{prop}
\label{FLMstableB}
Suppose that the fluid model is stable.  Then it is uniformly stable in the following two senses,
\begin{romannum}
\item
{\em Almost surely:}
 For $T\ge T_0$ and $\omega\in\Omega_0$,
\[
\lim_{\kappa\to\infty } \sup_{\|x\| = 1}   \| x^\kappa(T;x^\kappa,\omega) \|
=
0 .
\]
 \item
{\em In the $L_p$ sense:}
For $T\ge T_0$,
\begin{equation}
\lim_{\kappa\to\infty } \sup_{\|x\| = 1} \E[ \| x^\kappa(T;x^\kappa) \|^p]  = 0 .
\label{LpStable}
\end{equation}
\end{romannum}
\qed
\end{prop}

\begin{proof}[Proof of  \Proposition{Drift}]
Before proceeding it is helpful to review the Markov property:  Suppose that $\clC= \clC(X(0),X(1),\dots)$ is any random variable with finite mean.  We always have,
\[
\E_{X(n)}[\clC] = \E [\vartheta^n\clC \mid X(0),\dots, X(n)] =
\E [\vartheta^n\clC \mid  X(n)] ,
\]
where  $\clF_n\eqdef \sigma\{X(0),\dots,X(n)\}$,
$ n\ge 0$,  and
$\vartheta^n\clC$ denotes the random variable,
\[
\vartheta^n\clC= \clC(X(n),X(n+1),\dots)
\]
We apply the Markov property with $n=1$ and $\clC=\sum_{t=0}^{ \|X(0)\| T} \| X(t)\|^p $.  In this case we have,
\[
\vartheta^1\clC =
 \sum_{t=1}^{\|X(1)\|T} \|X(t)\|^p
\]
\notes{DG: I think it should be $\|X(1)\|T+1$. The proof is not affected by it\\
spm: no, I think it is correct.  This really is the defn of the shift.}
so that the Markov property gives,
\[
 V_p(X(1))
=
 \E_{X(1)}\bigl[  \clC \bigr]
=
 \E \Bigl[  \sum_{t=1}^{\|X(1)\|T} \|X(t)\|^p \mid \clF_1\Bigr].
\]

Applying the transition matrix to $V_p$ gives $PV_p(x)  =  \E_{x}[V_p(X(1))]$, so that
\[
 PV_p\, (x)  =
\E_{x}\Bigl[ \sum_{t=1}^{\|X(1)\|T} \|X(t)\|^p \Bigr]  ,\qquad
x\in\nat^N.
\]
The sum within the expectation on the right hand side is almost the same as used in the definition of $V_p$.  However, instead of summing from $0$ to $T\|X(0)\|$, we are summing from $1$ to $T\|X(1)\|$.  Consequently,  writing  $y=X(1)$, we have the expression,
\begin{equation}
PV_p(x)  =  V_p(x) - \|x\|^p +  \E_{x}\Bigl[ \sum_{t=\|x\|T+1}^{\|y\|T} \|X(t)\|^p \Bigr]  ,\qquad x\in\nat^N,
\label{PVfluidStability}
\end{equation}
where the sum is interpreted as negative when $\|x\| \ge \|y\|$.
Under the assumption that $\bfmL$ is a bounded sequence we obtain
the bound $PV_p \le V_p - \|x\|^p + b_p$,  where $b_p$ is the supremum over $x$ on the expectation on the right hand side of
\eqref{PVfluidStability}. We now argue that $b_p<\infty$. Recall that
under assumption (\ref{skip}) the increments of $\bfmX$ are bounded,  $| \|x\|T-\|X(1)\|T |\le b_0T$.  This combined with (\ref{LpStable}) implies that $b_p$ is indeed finite.
\notes{DG: I have added a sentence explaining why $b$ is finite.
\\
spm: great! I have made a few slight adjustments.}

To complete the proof we now establish \eqref{polyGrowth}.  For this we apply  \eqref{b0X}, which implies that $V_p$ satisfies the pair of bounds,
\[
\sum_{t=0}^{\|x\|T}  \|(x-b_0t)_+\|^p
\le
V_p(x)
\le
\sum_{t=0}^{\|x\|T}  \|x+b_0t\|^p,  \qquad   \ x\in\nat^N.
\]
This implies \eqref{polyGrowth}.
\end{proof}

\section{Scheduling policy}
\label{section:network}

The remainder of the paper is devoted to establishing the main result, Theorem~\ref{theorem:MainResult}.  Throughout the remainder of the paper we restrict attention to the KSRS model \eqref{RSQ} in continuous time.    The network is assumed symmetric with
$\lambda_1 =\lambda_3=1$,  $\mu_1=\mu_3=\infty$, and $\mu_2 =\mu_4$ finite.    The traffic intensity of each server is denoted,
\[
\begin{aligned}
 \rho_1 &\eqdef \lambda_1\mu_1^{-1}+\lambda_3\mu_4^{-1} = \mu_4^{-1}
\\
 \rho_2 &\eqdef \lambda_1\mu_2^{-1}+\lambda_3\mu_2^{-1} = \mu_2^{-1}
\end{aligned}
\]
To prove Theorem~\ref{theorem:MainResult} we construct a particular non-idling stationary policy.

The proofs of Lemmas~\ref{prop:psi} and \ref{lemma:parameters}  are contained in the  Appendix. Recall that the function $\Psi$ is defined in Theorem~\ref{theorem:MainResult}.
\begin{lemma}
\label{prop:psi}
For each constant $c>1$, the function $\Psi(c s)/\Psi(s)$ is a strictly increasing function. Moreover,
\begin{align*}
\sum_{1\leq m<\infty}m^2\Big[{\Psi(m^\eta)\over\Psi(c m^\eta)}\Big]^{1\over 4}<\infty,
\end{align*}
for any $\eta$ satisfying,
\begin{align}
\eta > \sqrt{\frac{12}{\log (c)}}
\label{eq:eta}
\end{align}
\end{lemma}

Throughout the remainder of the paper we fix $\mu_2,\mu_4$ so that the conditions specified in Lemma~\ref{lemma:parameters} are satisfied.
\begin{lemma}
\label{lemma:parameters}
Define $
\beta_1 =  {1\over 4}\gamma_{24} $ and $\beta_2 = 4\gamma_{24} $ with,
\[
  \gamma_2 = {\rho_2\over 1-\rho_2},  \quad
\gamma_4 = {\rho_1\over 1-\rho_1}, \quad
\gamma_{24} = \gamma_2\gamma_4,  \quad
 \gamma  = \gamma_4+\gamma_{24} .
\]
The parameters $\mu_2,\mu_4$ can be chosen so that $\rho_i\in (1/2,1)$ for $i=1,2$,  $\beta_2>\beta_1>1$,
and the condition (\ref{eq:eta}) is satisfied with $c=\beta_1$ and
$
\eta = \log( \beta_1)/ \log(\beta_2)$.
\end{lemma}

We define for $s,n\ge 1$,
\begin{equation}
\label{eq:Phi*}
\begin{aligned}
\Psi^*(s) &=(\Psi(\beta_1s^{\eta })/\Psi(s^{\eta }))^{1\over 4}
\\[.2cm]
 \psi(n)&=\Psi^*(n)/\Psi^*(n+1).
\end{aligned}
\end{equation}
By \Lemma{prop:psi}, $\Psi^*$ is a strictly increasing function, so
$\psi(n)\in (0,1)$.   The scheduling decisions are parametrized by
$\psi=\{\psi(n)\}$ and are made only at the sampling instances   $\{\tau_n, n\geq 0\}$ introduced in Section~\ref{stable}.

The scheduling decisions for server 1 at time $\tau_n$ are defined as follows:
\begin{romannum}
\item
If only one of the two buffers at server~1 contain jobs
 (buffer 1 or buffer~4), the server works on the first job in the non-empty buffer.  That is,  the policy is \textit{non-idling}.

\item
If both buffer 1 and buffer~4 are non-empty, and buffer~2 is also non-empty,  then buffer 1 receives strict priority at server 1.  Since $\mu_1=\infty$, all of the $Q_1(\tau_n)$ jobs are instantly sent to
buffer~2.

\item
If both buffer 1 and buffer~4 are non-empty, and buffer~2 is empty   then the scheduling decision depends on whether $\tau_n$ corresponds to an arrival into buffer 1 or not. If it does not correspond to an
arrival into buffer 1, then server 1 works on the first  job in
buffer~4.

\item
The policy is randomized if both buffer 1 and buffer~4 are non-empty, buffer~2 is empty, and  $\tau_n$ corresponds to an arrival into buffer 1.   In this case,  with probability $\psi(m)$ the server works on the first job in buffer~4,  and with probability $1-\psi(m)$ it works on the jobs in buffer 1.  This choice is made independently from any other randomness in the network.
In the second case, since the service rate is infinite, all of the  jobs in buffer 1 are instantly sent to buffer~2.
\end{romannum}
The scheduling decisions in server 2 are defined analogously.

With a slight abuse of notation we denote this scheduling policy by
$\psi$. We denote by  ${\cal Q}=(Q(0),\mu_2,\mu_4,\psi)$
the queueing network together with the scheduling policy $\psi$ and the initial state $Q(0)$.

The conclusions of \Proposition{prop:rates} are immediate from the discussion in Section~\ref{stable}.
\begin{prop}\label{prop:rates}
Under the scheduling policy $\psi$ the process $\bfmQ$ and the embedded process $\{X(n)\eqdef Q(\tau_n)\}$ are Markovian.
The chain $\bfmX$ satisfies the skip-free property \eqref{skip} as well as the irreducibility condition \eqref{irr}.  If $\bfmX$ possesses an invariant measure $\pi$ then it is necessarily unique, and satisfies  $\pi(x^*)>0$ with,
\begin{equation}
  x^*\eqdef (0,0,0,1)^\transpose.
\label{xstar}
\end{equation}
\end{prop}





The following  result is used to translate properties of  invariant measures to the process level properties in the network.

\begin{prop}\label{prop:FiniteExpectation}
Suppose the queueing network ${\cal Q}$ is such that an associated
invariant measure $\pi$ exists,
and $\pi(\Psi)$ is finite.     Then for every constant $\delta >0$.
\begin{align}\label{eq:FiniteExpectation}
\limsup_{s\rightarrow\infty}\Psi(\delta s)\pr_{x^*}({\|Q(s)\|\over
s}>\delta )<\infty.
\end{align}
Moreover, for every constant $\delta >0$
\begin{align}\label{eq:supratio}
\limsup_{s\rightarrow\infty}\Psi(.5\delta s)\pr_{x^*}(\sup_{t\geq s}{\|Q(t)\|\over t}>\delta )<\infty.
\end{align}
\end{prop}

\begin{proof}
By Proposition~\ref{prop:rates}, we have $\pi(x^*)>0$.
Applying Markov's inequality gives,
\begin{align*}
\Psi(\delta s)\pr_{x^*}(\|Q(s)\|>\delta s)&=\Psi(\delta s)\pr_{x^*}(\Psi(\|Q(s)\|)>\Psi(\delta s)) \\
&\leq \E_{x^*}[\Psi(\|Q(s)\|)]\\
&\stackrel{{\rm (*)}}{\leq} \pi^{-1}(x^*)\sum_{y}\E_y[\Psi(\|Q(s)\|)]\pi(y) \\
&= \pi^{-1}(x^*)\E_{\pi}[\Psi(\|Q(0)\|)],
\end{align*}
where in $(*)$ we simply use $1=\pi^{-1}(x^*)\pi(x^*)\leq \pi^{-1}(x^*)\sum_y\pi(y)$.
Then (\ref{eq:FiniteExpectation}) follows immediately.

We now establish (\ref{eq:supratio}).
Fix  $\delta >0$ and let $\hat \delta=\delta/4$.
Consider any $s>0$ and let  $s_k=(1+\hat\delta)^ks, k\geq 0.$ We use the following obvious identity,
\begin{align*}
\sup_{s_k\leq t\leq s_{k+1}}{\|Q(t)\|\over t}\leq {\|Q(s_k)\|\over
s_k}+{A_1(s_k,s_{k+1})+A_3(s_k,s_{k+1})\over s_k},
\end{align*}
with $A_i(s_k,s_{k+1})$ denoting the arrivals to buffer $i$ in the interval $(s_k,s_{k+1}]$.
Hence,
\begin{align*}
\pr(\sup_{t\geq s}{\|Q(t)\|\over t}>\delta )\leq \sum_{k\geq
0}\pr({\|Q(s_k)\|\over s_k}\geq \delta/2) +\sum_{k\geq
0}\pr({A_1(s_k,s_{k+1})+A_3(s_k,s_{k+1})\over s_k}\geq \delta/2).
\end{align*}
To complete the proof it suffices to show that
\begin{align}
&\limsup_{s\rightarrow\infty}\Psi(.5\delta s)\sum_{k\geq 0}\pr({\|Q(s_k)\|\over s_k}\geq \delta/2)<\infty, \label{eq:prbound} \\
&\limsup_{s\rightarrow\infty}\Psi(.5\delta s)\sum_{k\geq 0}\pr({A_1(s_k,s_{k+1})+A_3(s_k,s_{k+1})\over s_k}\geq \delta/2)
<\infty, \label{eq:arrivalbound}
\end{align}


Applying (\ref{eq:FiniteExpectation}), we have that
\begin{align*}
\sum_k\pr_{x^*}({\|Q(s_k)\|\over s_k}\geq \delta/2) =
O\Big(\sum_k[\Psi(.5\delta s_k)]^{-1}\Big) &=O\Big(\sum_k[\Psi(.5\delta
(1+\hat\delta)^ks)]^{-1}\Big)
\\
&=O\Big(\sum_k  \bigl(.5\delta (1+\hat\delta)^ks\bigr)^{-\log(.5\delta(1+\hat\delta)^ks)}\Big) \\
\end{align*}
Now for all $s$ such that $.5\delta s>1$ we have
\begin{align*}
\sum_k{1\over (.5\delta (1+\hat\delta)^ks)^{\log(.5\delta(1+\hat\delta)^ks)}}
&\le \sum_k{1\over (.5\delta s)^{\log(.5\delta(1+\hat\delta)^ks)}}\\
&=\sum_k{1\over (.5\delta s)^{\log(.5\delta s)+k\log(1+\hat\delta)}}\\
&=(.5\delta s)^{-\log(.5\delta s)}{1\over 1-(.5\delta s)^{-\log(1+\hat\delta)}}\\
&=\Psi^{-1}(.5\delta s){1\over 1-(.5\delta s)^{-\log(1+\hat\delta)}}.
\end{align*}
Since the limit of the second component of the product above is equal to unity, we obtain
that (\ref{eq:prbound}) holds.

We now establish (\ref{eq:arrivalbound}).
Using the large deviations bound given
by Lemma~\ref{lemma:LDexponential} in the Appendix, we have
\begin{align*}
\pr(A_1(s_k,s_{k+1})\geq (\delta /2)s_k)&=\pr(A_1(s_k,s_{k+1})-(s_{k+1}-s_k)\geq (\delta /2)s_k- (s_{k+1}-s_k)) \\
&\leq\pr(A_1(s_k,s_{k+1})- (s_{k+1}-s_k)\geq ({\delta \over 2\hat\delta}-1)(s_{k+1}-s_k))  \\
&=\pr(A_1(s_k,s_{k+1})- (s_{k+1}-s_k)\geq (s_{k+1}-s_k))  \\
&\leq \exp(-\Omega(s_{k+1}-s_k)) \\
&\leq \exp(-\Omega(\hat\delta(1+\hat\delta)^k s)) \\
\end{align*}
We also have
\begin{align*}
\sum_{k\geq 0}\exp(-\Omega(\hat\delta(1+\hat\delta)^k s))\leq\exp(-\Omega(s)).
\end{align*}
Since $\Psi$  grows slower than exponentially we conclude that (\ref{eq:arrivalbound}) holds.
\end{proof}

\section{Stability of the fluid limit model}
\label{section:stability}

The goal of this section is to establish stability of the fluid limit model obtained from the scheduling policy $\psi$.
Specifically, we establish the following result which is Part 2 of Theorem~\ref{theorem:MainResult}.

\begin{theorem}\label{theorem:stability1}
There exists a constant $c_0>0$ such that for every $\epsilon>0$
\begin{align*}
\limsup_{\|x\|\rightarrow\infty}\pr_x\Big({\|Q(c_0\|x\|)\|\over
\|x\|}>\epsilon\Bigr)=0.
\end{align*}
\end{theorem}

\begin{proof}
The result is obtained by combining Proposition~\ref{prop:Q4} with Lemma~\ref{lemma:Q4} that follow.
\end{proof}

The proofs in this section rely on regeneration arguments, based
on the stopping times defining the first emptying times for the four buffers in the network,
\begin{equation}
T_i=\inf\{t\geq 0:Q_i(t)=0 \}
\quad
\barT_i=\inf\{t\geq 0:Q_j(t)=0\ \hbox{for all \ $j\neq i$, $i=1,\dots,4$}\} .
\label{stop}
\end{equation}
The following feature of the model is crucial to obtain regenerative behavior:
 Although the policy is not a priority policy, its behavior at buffers~4 and 2 is similar.
 For example, since the policy is non-idling and the service rate at buffer 1 is
 infinite we can conclude that buffer~4 receives full capacity while this buffer is non-empty.
\notes{spm 12/18: minor revisions here}

\Lemma{mm1Like} follows from these observations, and the specification of the policy that gives priority to  jobs in
buffer~3 while buffer~4 is non-empty.
\begin{lemma}
\label{mm1Like} Under the policy $\psi$ with $Q_4(0)\ge 1$,  the
contents of buffer~4 evolves  as an M/M/1 queueing system with
arrival rate $1$ and service rate $\mu_4$, until the first time $Q_4$ becomes zero.
\end{lemma}

Proposition~\ref{prop:Q4} concerns  the special case in which all the jobs are initially in buffer~4.

\begin{prop}\label{prop:Q4}
Let $x^n=(0,0,0,n)=nx^*$ and $\tau=1/(\mu_4-1)$. For every $\epsilon>0$
\begin{align*}
\lim_n\pr_{x^n}\Big(\|Q(\tau n)\|\leq \epsilon n\Big)=1.
\end{align*}
\end{prop}

\begin{proof}
We begin with an application of \Lemma{mm1Like} to conclude that
$Q_4(t)$ corresponds with an M/M/1 queueing system up to time  $T_4$.     Applying Lemma~\ref{lemma:emptyingTimeMM1}
we obtain for every $\epsilon>0$
\begin{align}\label{eq:tau4}
\lim_n\pr_{x^n}\Big(|T_4-{n\over \mu_4-1}|\leq \epsilon n\Big)=1.
\end{align}

We now analyze the queue-length processes in buffers~1 and 2 during $[0,T_4)$.  For each $k\ge 1$, with probability $(1-\psi(k))\prod_{1\leq i\leq k-1}\psi(i)$,  priority is given to jobs in buffer~4 up until the time of the $k$-th arrival.   At this time the $k$ jobs at buffer~1 are immediately transfered to buffer~2.  While buffer~2 is non-empty, full priority is given to jobs in buffer 1 over jobs in buffer~4, so that any new arrivals are  sent to buffer~2 instantaneously.
At the first instance $\tau_2(1)<T_4$ that the buffer~2 empties
(assuming one exists), this process repeats.


Hence the policy $\psi$ induces the following behavior:  there is an alternating  sequence $R_1,L_1,R_2,L_2,\ldots$, where $R_l$
corresponds to time-intervals before $T_4$ when
buffer~2 is empty, and $L_l$ corresponds to time-intervals before $T_4$ when buffer~2 is non-empty.   The sequences $R_l, l\geq 1$ and $L_l\geq l$ are i.i.d., and the length of $R_l$ is independent from $R_{l'},L_{l'}, l'\leq l-1$.

Let $\tau_l=\sum_{j\leq l}(R_j+L_j)$.  The queue length in buffer 1 at the end of the time period corresponding to $R_l$ (that is
at time $\tau_{l-1}+R_l$) is equal to $k$ with
probability $(1-\psi(k))\prod_{1\leq i\leq k-1}\psi(i)$. Let $X_l$ represent this queue length
$Q_1(t)$ at time $t=\tau_{l-1}+R_l$.
 We have
\begin{align*}
\E[X_l^2]&=\sum_{k\ge 1}k^2(1-\psi(k))\prod_{1\leq i\leq k-1}\psi(i)\\
&<\sum_{k\ge 1}k^2\prod_{1\leq i\leq k-1}\psi(i)\\
&=\sum_{k\ge 1}{k^2\Psi^*(1)\over \Psi^*(k)}
\end{align*}
Applying the second part of \Lemma{prop:psi}, this sum is finite, namely $X_l$ has a finite
second moment. Then $R_l$ has finite second moment as well since, conditioned on $X_l=k$,  it is a sum of $k$ i.i.d. random variables with $\Exp(\lambda)$ distribution.
Conditioning on $X_l=k$, the length of $L_{l}$ represents the time to empty an M/M/1 queueing system with $k$ initial jobs, arrival rate $1$,
and service rate $\mu_2$. Hence  $\E[L_l^2|X_l=k]=O(k^2)$ by
Lemma~\ref{lemma:emptyingTimeMM1}. Since $X_l$ has a finite second moment, so does $L_l$. We conclude that $\tau_l-\tau_{l-1}=R_l+L_l$ has a finite second moment.

Assume now that $T_4=\infty$ by placing infinitely many jobs in buffer~4.  Given any positive $t>0$, let $l^*(t)$ be the unique index such that $t\in [\tau_{l^*(t)-1},\tau_{l^*(t)}]$.
Applying Smith's Theorem for regenerative processes (see Theorem 3.7.1. in Resnick~\cite{ResnickStochasticProcesses}),
for every $m\geq 0$
\begin{align*}
\lim_{t\rightarrow\infty}\pr(Q_1(t)+Q_2(t)\geq m)={\E[\int_0^{\tau_1}1\{Q_1(t)+Q_2(t)\geq m\}dt]\over \E[\tau_1]},
\end{align*}
where $\tau_1=R_1+L_1$.
Applying the Cauchy-Schwartz inequality we obtain,
\begin{align*}
\E[\int_0^{\tau_1}1\{Q_1(t)+Q_2(t)\geq m\}dt&\leq \E[1\{\sup_{0\leq t\leq \tau_1}Q_1(t)+Q_2(t)\geq m\}\tau_1] \\
&
\le \sqrt{\pr(\sup_{0\leq t\leq \tau_1}Q_1(t)+Q_2(t)\geq m)\E[\tau_1^2] } \\
&
\leq \sqrt{\pr(A_1(\tau_1)\geq m)\E[\tau_1^2]},
\end{align*}
Observe that $A_1(\tau_1)=X_1+A(R_1,\tau_1)$.
Conditioning on $X_1=k$ and applying Lemma~\ref{lemma:emptyingTimeMM1} gives $\E[A(R_1,\tau_1)]=O(k)$, and since the mean of $X_1$ is finite,  $\E[X_1]<\infty$, we obtain $\E[A_1(\tau_1)]<\infty$.

Markov's inequality, gives the bound $\pr(A_1(\tau_1)\geq m)=O(1/m)$, and since $\E[\tau_1^2]<\infty$ we conclude,
\begin{align*}
\lim_{t\rightarrow\infty}\pr_{x^\infty}(Q_1(t)+Q_2(t)\geq m)=O({1\over m}).
\end{align*}
Now we recall (\ref{eq:tau4}) and use independence of interarrival and service times to conclude
that
\begin{align*}
\lim_n\pr_{x^n}(Q_1(T_4)+Q_2(T_4)\geq m)=O({1\over m}) .
\end{align*}
Since $Q_3(T_4)+Q_4(T_4)=0$, we obtain an apparently weaker result that for every $\epsilon>0$
\begin{align*}
\lim_n\pr_{x^n} \bigl ( \|Q(T_4)\| \leq \epsilon n \bigr)=1.
\end{align*}
It remains to relate $\|Q(T_4)\|$ to $\|Q({n\over \mu_4-1})\|$. We
use the fact that the difference between $\|Q({n\over \mu_4-1})\|$
and $\|Q(T_4)\|$ is at most the total number of arrivals and
departures during the time interval between ${n\over \mu_4-1}$ and
$T_4$. Using the large deviations bound
Lemma~\ref{lemma:LDexponential} applied to arrival processes to
buffers 1 and 3 and service processes in buffers 2 and 4, and
combining with (\ref{eq:tau4}), we obtain
\begin{align*}
\lim_n\pr_{x^n} \bigl( \bigl\|Q \bigl({n\over \mu_4-1}  \bigr ) \bigr \| \leq \epsilon n\bigr)=1.
\end{align*}
\end{proof}

We now assume that $x$ is not of the form $(0,0,0,n)$ and
complete the proof of Theorem~\ref{theorem:stability1} by establishing the following lemma.

Recall that $\{\barT_i\}$ are defined in \eqref{stop}, given
$Q(0)=x$. We let $T$ denote the minimum,
\[
T=\min(\barT_2,\barT_4).
\]
\begin{lemma}
\label{lemma:Q4} For some constants $c_1,c_2>0$,
\begin{align*}
\liminf_{\|x\|\rightarrow\infty}\pr_x\big(T \le   c_1\|x\| \
\hbox{\it and}\ \|Q(T)\| \le c_2 \|x\| \big)=1.
\end{align*}
\end{lemma}

\begin{proof} 
To obtain bounds on the probability in the lemma  we claim it is enough to consider the special case in which one of the servers has no jobs at time $t=0$.

Suppose both of the servers are initially non-empty, let $x=Q(0)$,
and consider the following cases.  If $Q_4(0)=0$ then this violates
the right-continuity assumption since we then have $Q_1(0) >0$,  and
all of the jobs in buffer 1 proceed immediately to buffer~2. In the
second case   $Q_4(0)>0$, and we can apply \Lemma{mm1Like} to deduce
that $T_4(x) = O(\|x\|)$ with probability approaching unity for
large $\|x\|$.  At time $T_4$ buffer~3 is empty and all of the jobs
in buffer 1 immediately proceed to buffer~2.   Once again we arrive
at a state $Q(T_4)$ corresponding to an empty server.

In the remainder of the proof we assume that one server is empty at time $t=0$;  Without loss of generality this is server~2, so that  $x_2+x_3=0$.  Again, if in addition $x_4=0$ then only buffer 1 is
non-empty, and all the jobs in buffer 1 are sent to buffer~2, giving $T=0$.

Otherwise,  suppose that  $x_2+x_3=0$, $x_4>0$, and consider the
queue length process at buffer~4. Applying \Lemma{mm1Like} we
conclude that buffer~4 evolves as  an M/M/1 queueing system, and
Lemma~\ref{lemma:emptyingTimeMM1} implies that the emptying time
$T_4$ for buffer~4 is less than $c_1\|x\|$, with probability approaching unity as
$n\rightarrow\infty$, for some constant $c_1$. At $t=T_4$ buffer~3
is empty as well and all the jobs in buffer 1 (if any) instantly
proceed to buffer~2. Hence buffers 1,3 and 4 are all empty at time $
T_4$, so that $T=T_4$.

This gives a uniform bound on the probability $ \pr_x\big(T \le
c_1\|x\|\bigr)$ for initial conditions corresponding to one empty
server. To obtain a uniform bound on     $ \pr_x\big( \{T \le
c_1\|x\|\}\cap \{ \|Q(T)\| \le c_2 \|x\| \} \bigr)  $ for some $c_2$
we apply the large deviations bound of
Lemma~\ref{lemma:LDexponential} to the arrival processes along with
the bound $\|Q(t)\|\le \|x\| + \|A(t)\|$.
\end{proof}

\section{Lower bounds on the tail probability}
\label{section:instability}

The goal of this section is to prove Part~3 of
Theorem~\ref{theorem:MainResult}.   In light of Theorem~\ref{geotailRelax} it suffices to establish (\ref{eq:expmomentInfinite}),
and then the lack of exponential ergodicity will follow.
\notes{*.*}

Recall the definition of the first  emptying time $T_4$ from
\eqref{stop}.   We fix a constant $\epsilon>0$,  and define the
following events given $Q(0)=x$,
\begin{align*}
{\cal E}_4&=1\Big\{{(1-\epsilon)x_4\over \mu_4-1}\leq T_4\leq {(1+\epsilon)x_4\over \mu_4-1}\Big\}, \\
{\cal E}_1&=1\Big\{{(1-\epsilon)^2  x_4\over \mu_4-1}\leq
Q_2(T_4)\leq {(1+\epsilon)^2  x_4\over \mu_4-1} \wedge
Q_i(T_4)=0,~i\neq 2\Big\}.
\end{align*}

Like \Lemma{lemma:Q4}, the following result compares an emptying time for the
stochastic model with the emptying time for a fluid model.
\begin{lemma}\label{lemma:arrivalst_m}
Consider a non-zero initial state $x=(x_1,x_2,x_3,x_4)$ satisfying
$x_i=0,i\neq 4$. Then for any $0<\epsilon<1/10$
\begin{align}
\pr_x({\cal E}_4 \cap {\cal E}_1)\geq \Theta((\Psi^*(2\gamma_4
x_4))^{-1}). \label{eq:Q1}
\end{align}
\end{lemma}

\begin{proof}
Consider first a modified scheduling policy $\tilde\psi$ that always
gives  priority to jobs in buffer~4 in server 1 and buffer~3 in
server 2.  For the   policy $\tilde\psi$ we let
$\{\tilde T_i \}$  denote the draining times for the four buffers.
For this policy all  jobs in buffer~3 are transfered  to buffer~4
instantaneously so that buffer 4 operates as an M/M/1 queueing
system for all $t\ge 0$.   Applying the bound (\ref{eq:LDonT}) from
Lemma~\ref{lemma:emptyingTimeMM1}, the stopping time $\tilde T_4$
satisfies,
\begin{align}
\pr\Big({(1-\epsilon)x_4\over \mu_4-1}\leq \tilde T_4\leq
{(1+\epsilon)x_4\over \mu_4-1}\Big) \geq 1-\exp(-\Omega(x_4)).
\label{eq:tildeT4x}
\end{align}
Let $m_1(x)$ denote the number of arrivals to buffer 1 before buffer~4 becomes empty for $\tilde\psi$.  That is,
\begin{align*}
m_1(x)=A_1(\tilde T_4).
\end{align*}
Using Lemma~\ref{lemma:LDexponential} applied to the arrival process, combined with (\ref{eq:tildeT4x}), we obtain
\begin{align}
\pr\Big({(1-\epsilon)^2  x_4\over \mu_4-1}\leq m_1(x) \leq
{(1+\epsilon)^2  x_4\over \mu_4-1}\Big)\geq 1-\exp(-\Omega(x_4)).
\label{eq:m1x}
\end{align}

Now we return to the original scheduling policy $\psi$.  For every
$m$, conditioned on $m_1(x)=m$, the  probability that at every
arrival instance the priority was given to jobs in buffer~4 is
$\prod_{1\leq i\leq m}\psi(i)=(\Psi^*(m))^{-1}$.   If it is indeed the case
that at every arrival into buffer 1 the priority was given to
buffer~4 for all arrivals up to the  $m_1(x)$-th, then $T_4=\tilde
T_4$ and $Q_1(T_4)=m_1(x)$. If in addition $m=m_1(x)$ satisfies the
bound
\begin{align*}
m\leq {(1+\epsilon)^2  x_4\over \mu_4-1}={(1+\epsilon)^2\rho_1x_4\over 1-\rho_1}
<2\gamma_4 x_4,
\end{align*}
then, by monotonicity of $\Psi^*$, we have $(\Psi^*(m))^{-1}\geq
(\Psi^*(2\gamma_4 x_4))^{-1}$. If during the time interval $[0,T_4]$
server 1 is processing only jobs in buffer~4, then at time $T_4$
server 2 is empty. Since at time $T_4$ buffer~4 also becomes empty,
all the $Q_1(T_4)$ jobs in buffer 1  instantly arrive into buffer~2
at time $T_4$, and all the other buffers become empty at $T_4$.
\notes{spm: this is said so often - should we add another lemma??  - I haven't done so, but we can think about it in the revision.} That
is $Q_i(T_4)=0, i\neq 2$ and $Q_2(T_4)=m_1(x)$. These arguments
combined with  (\ref{eq:tildeT4x}) and (\ref{eq:m1x}) give,
\begin{align*}
\pr_x({\cal E}_4 \cap {\cal E}_1)=(\Psi^*(2\gamma_4
x_4))^{-1}-2\exp(-\Omega(x_4))=\Theta((\Psi^*(2\gamma_4 x_4))^{-1}),
\end{align*}
where the last equality follows from subexponential decay of $\Psi^*$.
This completes the proof.
\end{proof}

Define the stopping time,
\begin{equation}
T_{24}=\inf\{t>T_4:Q_2(t)=0\},
\label{T42}
\end{equation}
and for a given $x=Q(0)$ consider the events,
\begin{align*}
{\cal E}_{42}&=1\Big\{{(1-\epsilon)^3  x_4\over (\mu_4-1)(\mu_2-1)}\leq T_{24}-T_4\leq
{(1+\epsilon)^3  x_4\over (\mu_4-1)(\mu_2-1)}\Big\}, \\
{\cal E}_{13}&=1\Big\{{(1-\epsilon)^4  x_4\over (\mu_4-1)(\mu_2-1)}\leq Q_4(T_{24})\leq
{(1+\epsilon)^4  x_4\over (\mu_4-1)(\mu_2-1)} \wedge Q_i(T_{24})=0, ~i\neq 4\Big\}.
\end{align*}

\begin{lemma}\label{lemma:arrivalst_m'}
Consider a starting state $x=(x_1,x_2,x_3,x_4)$ satisfying
$x_i=0,i\neq 4$. Then for every $\epsilon>0$
\begin{align}
\pr_x({\cal E}_{42}\cap {\cal E}_{13}\cap {\cal E}_4\cap {\cal
E}_1)= \Theta((\Psi^*(4\gamma_{24} x_4))^{-2}). \label{eq:Q3}
\end{align}
\end{lemma}

\begin{proof}
The event ${\cal E}_1$ is equivalent to  $\{ Q(T_4)=(0,z,0,0)
\hbox{  for some } z\in \clS_4\}$, where $\clS_4$ is the set
of integers,
\begin{align*}
\clS_4=\{z\in\nat:{(1-\epsilon)^2  x_4\over \mu_4-1}\leq z\leq {(1+\epsilon)^2  x_4\over \mu_4-1}\}.
\end{align*}
Consequently,
\begin{align*}
\pr\Big({\cal E}_{42}\cap {\cal E}_{13}\cap {\cal E}_4\cap {\cal
E}_1\Big)
&=\sum_{z\in \clS_4}\pr\Big({\cal E}_{42}\cap {\cal E}_{13}\cap (Q(T_4)=(0,z,0,0)) \cap {\cal E}_4\Big) \\
&=\sum_{z\in \clS_4}\pr\Big({\cal E}_{42}\cap {\cal E}_{13}\big|
Q(T_4)=(0,z,0,0)\Big) \pr\Big(Q(T_4)=(0,z,0,0) \cap {\cal E}_4\Big),
\end{align*}
where in the second equality we use the Markovian property of the scheduling policy $\psi$. Applying Lemma~\ref{lemma:arrivalst_m} but
interchanging buffer~4 with buffer~2 and buffer~3 with buffer 1, we obtain  for every $z$.
\begin{align*}
\pr\Big({\cal E}_{42}\cap {\cal E}_{13}\big|  Q(T_4)=
(0,z,0,0)\Big)\geq \Theta((\Psi^*(2\gamma_2 z))^{-1}).
\end{align*}
For every $z\in \clS_4$, given the bound $\epsilon<1/10$, we have $z\leq 2\gamma_4x_4 $. By monotonicity of
$\Psi^*$ and the definition $\gamma_{24} = \gamma_2\gamma_4$
we obtain $[\Psi^*(2\gamma_2 z) ]^{-1}\geq [\Psi^*(4\gamma_{24} x_4)]^{-1}$.

Finally, we note that,
\begin{align*}
\sum_{z\in \clS_4} \pr\Big(  \{Q(T_4)=(0,z,0,0) \}\cap {\cal
E}_4\Big) = \pr\Big({\cal E}_1\cap {\cal E}_4\Big)
\end{align*}
which is at least $\Theta((\Psi^*(2\gamma_4 x_4))^{-1})$ by Lemma~\ref{lemma:arrivalst_m}.
Now $\rho_2>1/2$ implies $\gamma_2>1$ which gives
$\gamma_{24}=\gamma_2\gamma_4>\gamma_4$. Combining these bounds we obtain the desired bound (\ref{eq:Q3}).
\end{proof}

\begin{lemma}\label{lemma:fullcycle}
Given a state $x$ such that $x_i=0,i\neq 4$, there exists a random time $T$ such that,
\begin{align*}
\pr\Big({\gamma x_4\over 2}\leq  T\leq & 2\gamma x_4~ \wedge~ {
\gamma_{24}\over 2\gamma }\leq {Q_4(T)\over T}\leq {2
\gamma_{24}\over\gamma } ~ \wedge~ Q_i(T)=0, i\neq 4\Big) \breakMed
& \geq \Theta((\Psi^*(4\gamma_{24}x_4))^{-2}).
\end{align*}
\end{lemma}

\begin{proof}
The event ${\cal E}_4\cap{\cal E}_{42}$ implies
\begin{align*}
((1-\epsilon)\gamma_4+(1-\epsilon)^3\gamma_{24})x_4\leq T_{24}
\leq ((1+\epsilon)\gamma_4+(1+\epsilon)^3\gamma_{24})x_4,
\end{align*}
  the event ${\cal E}_{13}$ implies
\begin{align*}
(1-\epsilon)^4\gamma_{24}\leq Q_4(T_{24})\leq (1+\epsilon)^4,\gamma_{24}
\end{align*}
and the bound $\epsilon<1/10$ implies that
$1/2<{(1-\epsilon)^4\over (1+\epsilon)^4}<{(1+\epsilon)^4\over (1-\epsilon)^4}<2$.   The result is obtained from Lemma~\ref{lemma:arrivalst_m'} and using
$\gamma\eqdef \gamma_4+\gamma_{24}$.
\end{proof}

We now use Lemma~\ref{lemma:fullcycle} is used to obtain the following bound.

\begin{prop}\label{prop:manycycles}
There exist  constants $0<\alpha<1,~\delta>0$ such that for every $n\in\nat$
\begin{align}\label{eq:eq1}
\pr_{x^*}\Big[\sup_{t\geq \beta_1^n}{\|Q(t)\|\over
t}>\delta\Big]\geq \alpha^n\prod_{1\leq m\leq
n}(\Psi^*(\beta_2^{m+1}))^{-2}.
\end{align}
\end{prop}

\begin{proof}
For every $n\geq 1$ denote by ${\cal E}_n$ the event that the event described in Lemma~\ref{lemma:fullcycle}
occurs $n$ times in succession starting from the state $x=x^*$.
For each $m$ let $\sigma_m$ is the length of the time interval corresponding to the $m$-th event, and let $S_n= \sum_{m=1}^n  \sigma_m$.     The random time $S_n$ is only defined on the event ${\cal E}_n$.

Using Lemma~\ref{lemma:fullcycle}, the event ${\cal E}_n$ implies  that ${\gamma \over 2 }\leq \sigma_1\leq 2\gamma $,
and for each $2\leq m\leq n$,
\begin{equation}
\begin{aligned}
{\gamma \over 2 }Q_4(\sum_{j\leq m-1}\sigma_j)&\leq \sigma_{m}\leq 2\gamma Q_4(\sum_{j\leq m-1}\sigma_j),
\breakMed
{  \gamma_{24}\over 2\gamma }\sigma_{m}&\leq Q_4(\sum_{j\leq m}\sigma_j)\leq {2  \gamma_{24}\over \gamma }\sigma_{m}.
\end{aligned}
\label{eq:Tlower}
\end{equation}
Combining, we obtain
\begin{align}
\sigma_n\geq {\gamma_{24}\over 4}\sigma_{n-1}\geq {\gamma_{24}^2\over 4^2}\sigma_{n-2}\geq\cdots\geq  {\gamma_{24}^{n-1}\over 4^{n-1}}\sigma_1
\geq{\gamma_{24}^{n-1}\over 4^{n-1}}{\gamma \over 2 } ,
 \label{eq:T_m}
\end{align}
and
\begin{align}
\sigma_n\leq 4\gamma_{24}\sigma_{n-1}\leq (4\gamma_{24})^2 \sigma_{n-2}\leq\cdots\leq  (4\gamma_{24})^{n-1}\sigma_1
\leq (4\gamma_{24})^{n-1}2\gamma ,
\label{eq:T_mupper}
\end{align}
The chain of inequalities (\ref{eq:T_m}) implies  in particular that
\begin{align}
S_n \geq \sigma_n\geq {\gamma_{24}^{n-1}\over 4^{n-1}}{\gamma \over2 }. \label{eq:TTn}
\end{align}
Since $\beta_1=\gamma_{24}/4>1$ by Lemma~\ref{lemma:parameters},    the same chain of inequalities implies
\begin{align}
S_n =&\sum_m \sigma_m\leq \sigma_n\sum_{m\leq n} {\gamma_{24}^{-(n-m)}\over 4^{-(n-m)}}
\leq \sigma_n\sum_{j=1}^{\infty} {\gamma_{24}^{-j}\over 4^{-j}}
={\sigma_n\over 1-{4\over \gamma_{24}}}.
\label{eq:TTn1}
\end{align}
Also the event ${\cal E}_n$ implies
\begin{align*}
{\|Q(S_n)\|\over \sigma_n}\geq { \gamma_{24} \over 2\gamma },
\end{align*}
Then from (\ref{eq:TTn1}) we obtain
\begin{align*}
{\|Q(S_n)\|\over S_n}\geq { \gamma_{24} \over 2\gamma }(1-{4\over
\gamma_{24}}).
\end{align*}
Recalling the bound (\ref{eq:TTn}), we conclude that the event ${\cal E}_n$ implies
\begin{align*}
\sup \{ t^{-1}\|Q(t)\| :  t\geq {\gamma \over 2}({\gamma_{24}\over
4})^{n-1} \} \geq {\|Q(S_n)\|\over S_n} \geq { \gamma_{24} \over
2\gamma }(1-{4\over \gamma_{24}}).
\end{align*}

We have $\beta_1=\gamma_{24}/4<\gamma/2$
and conclude that the event
${\cal E}_n$ implies
\begin{align*}
\sup_{t\geq \beta_1^n}{\|Q(t)\|\over t}>\delta \eqdef{ \gamma_{24}
\over 2\gamma }(1-{4\over \gamma_{24}}).
\end{align*}
It remains to obtain the required lower bound on  $\pr_{x^*}[{\cal
E}_n]$.  We first obtain a bound on $\pr_{x^*}[{\cal E}_m|{\cal
E}_{m-1}], m=2,3,\ldots,n$. For convenience we set ${\cal
E}_0=1\{Q(0)=x^*\}$. For each $1\leq m\leq n$ the event ${\cal
E}_{m-1}$ via (\ref{eq:T_mupper}) and (\ref{eq:Tlower}) implies
\begin{align*}
Q_4(\sum_{j \leq m-1}\sigma_j)\leq {2 \over \gamma }\sigma_{m}
&\leq {2 \over \gamma }(4\gamma_{24})^{m-1}2\gamma =4^m\gamma_{24}^{m-1}.
\end{align*}
Let $\alpha>0$ be a constant hidden in the $\Theta(\cdot)$ notation in Lemma~\ref{lemma:fullcycle}.
Then, from Lemma~\ref{lemma:fullcycle}, we obtain
\begin{align*}
\pr_{x^*}[{\cal E}_m|{\cal E}_{m-1}]\geq
\alpha(\Psi^*((4\gamma_{24})4^m\gamma_{24}^{m-1}))^{-2}\geq \alpha(\Psi^*((4\gamma_{24})^{m+1}))^{-2}.
\end{align*}
Since this holds for every $n\geq m\geq  1$ and ${\cal E}_1\supset{\cal E}_2\supset\cdots\supset{\cal E}_n$ we obtain
\begin{align*}
\pr_{x^*}[{\cal E}_n]\geq \alpha^n\prod_{1\leq m\leq n}(\Psi^*((4\gamma_{24})^{m+1}))^{-2}.
\end{align*}
This concludes the proof of the proposition by substituting $\beta_2$ for $4\gamma_{24}$.
\end{proof}

We are ready to conclude the proof of the final part of Theorem~\ref{theorem:MainResult}.

\begin{proof}[Proof of Part 3 of Theorem~\ref{theorem:MainResult}]
Using the expression (\ref{eq:Phi*}) for $\Psi^*$ we have
\begin{align*}
\prod_{1\leq m\leq n}\Big(\Psi^*(\beta_2^{m+1})\Big)^2&=
\prod_{1\leq m\leq n}\Big({\Psi(\beta_1\beta_2^{\eta(m+1)})\over \Psi(\beta_2^{\eta(m+1)})}\Big)^{1\over 2} \\
&=\prod_{1\leq m\leq n}\Big({\Psi(\beta_1^{m+2})\over \Psi(\beta_1^{m+1})}\Big)^{1\over 2}\\
&<\Psi^{1\over 2}(\beta_1^{n+2}).
\end{align*}
>From Proposition~\ref{prop:manycycles}, we obtain that for all $n$,
\begin{align*}
(1/\alpha)^n\Psi^{1\over 2}(\beta_1^{n+2})\pr_{x^*}\Big[\sup_{t\geq
\beta_1^n}{\|Q(t)\|\over t}>\delta\Big] \ge 1.
\end{align*}
Using $\beta_1^n=s$ and finding a constant $\eta_2>0$ such that $(1/\alpha)=\beta_1^{\eta_2}$, we
obtain
\begin{align*}
\liminf_{s\rightarrow\infty}s^{\eta_2}\Psi^{1\over
2}(\beta_1^2s)\pr_{x^*}\Big[\sup_{t\geq s}{\|Q(t)\|\over
t}>\delta\Big] \geq 1.
\end{align*}
Now
\begin{align*}
s^{\eta_2}\Psi^{1\over 2}(\beta_1^2s)=s^{\eta_2}(\beta_1^2s)^{.5\log(\beta_1^2s)}
=\exp(.5\log^2 s+2\log\beta_1\log s+2\log^2\beta_1+\eta_2\log s).
\end{align*}
Observe, that the righ-hand side, as a function of $s$ is
\begin{align*}
o(e^{\log^2(.5\delta s)})=o(\Psi(.5\delta s))
\end{align*}
since the leading term in the first term is $e^{.5 \log^2 s}$ and in the second
term is $e^{\log^2 s}$.
We conclude
\begin{align*}
\liminf_{s\rightarrow\infty}\Psi(.5\delta s)\pr_{x^*}\Big[\sup_{t\geq s}{\|Q(t)\|\over t}>\delta\Big]=\infty.
\end{align*}
Applying the second part of Proposition~\ref{prop:FiniteExpectation} we obtain (\ref{eq:expmomentInfinite}).
The proof of Theorem~\ref{theorem:MainResult} is complete.
\end{proof}

\appendix

\section{Appendix}
In this section we provide proofs of some of the elementary results we have used above. Some of the
results are well known.

\begin{proof}[Proof of \Lemma{prop:psi}]
We have
\begin{align*}
{\Psi(c  s)\over \Psi(s)}={(c  s)^{\log (c  s)}\over s^{\log s}}=c  ^{\log (c  s)}s^{\log c  },
\end{align*}
which is a strictly increasing function for every $c  >1$. Now suppose $\eta>0$ is such
that (\ref{eq:eta}) is satisfied. Then
\begin{align*}
\sum_{1\leq m<\infty}m^2\Big[{\Psi(m^\eta)\over\Psi(c  m^\eta)}\Big]^{1\over 4}&=
\sum_{1\leq m<\infty}m^2{m^{{1\over 4}\eta\log m^{\eta}}\over (c  m)^{{1\over 4}\eta\log (c  m)^{\eta}}} \\
&=\sum_{1\leq m<\infty}{m^2\over c  ^{{1\over 4}\eta^2\log c  m}m^{{1\over 4}\eta^2\log c  }}.
\end{align*}
The condition (\ref{eq:eta}) implies ${1\over 4}\eta^2\log c -2>1$, which implies the result.

\end{proof}

\begin{proof}[Proof of Lemma~\ref{lemma:parameters}]
We fix a small $\delta>0$ and let $\mu_2=\mu_4=1+\delta$. Then
$\rho_2=\rho_1=1/(1+\delta),\gamma_2=\gamma_4=1/\delta,\gamma_{24}=1/\delta^2,
\beta_2=4/\delta^2,\beta_1=1/(4\delta^2),\eta={\log(1/\delta^2)+4\over\log(1/\delta^2)-4},
c=\beta_1=1/(4\delta^2).$ Then
\begin{align*}
\eta^2\log c=\Big({\log(1/\delta^2)+4\over\log(1/\delta^2)-4}\Big)^2\log ({1\over 4\delta^2}).
\end{align*}
As $\delta\rightarrow 0$, $\eta\rightarrow 1, c\rightarrow \infty$ and the condition (\ref{eq:eta})
is satisfied for sufficiently small $\delta>0.$ All the other constraints are satisfied trivially.

\end{proof}

We  use in this paper very crude large deviations type bounds. Of course in most
cases far more refined large deviations estimates are available \cite{large_deviations},
but those are not required for our purposes.

\begin{lemma}\label{lemma:LDexponential}
Let $N(t)$ be a Poisson process with parameter $\nu>0$. Then for every constant $\epsilon>0$
\begin{align*}
\pr(|N(t)-\nu t|>\epsilon t)=\exp(-\Omega(t)).
\end{align*}
\end{lemma}
Here the constants hidden in $\Omega$ include $\epsilon$.


\begin{lemma}\label{lemma:emptyingTimeMM1}
Consider an M/M/1 queueing system with parameters
$\lambda,\mu,\rho=\lambda/\mu<1$. Let $Q(t)$ denote the queue length
at time $t$ and let $T=\inf\{t:Q(t)=0\}$. Then $\E[T|Q(0)=n]=O(n)$
and $\E[T^2|Q(0)=n]=O(n^2)$. Moreover
\begin{align}\label{eq:LDonT}
\pr({T\over n\mu(1-\rho)} \in (1-\epsilon,1+\epsilon))=1-\exp(-\Omega(n)).
\end{align}
\end{lemma}

\begin{proof}
These are well known results from queueing theory. One quick way to
establish them is to observe that $T=\sum_{1\leq j\leq n}T_j$, where
$T_j$ is the first passage time from $j$ to $j-1$. That is
$T_j=\inf\{t:Q(t)=j-1|Q(0)=j\}$. The sequence $T_j$ is i.i.d. with
the distribution equal to the distribution of the busy period, which
is known to satisfy $\E[\exp(s T_j)]<\infty$ for some $s>0$, and has
the first moment  $1/(\mu(1-\rho))$. We immediately obtain
$\E[T|Q(0)=n]=n/(\mu(1-\rho))$ and $\E[T^2|Q(0)=n]=O(n^2)$ (actual
value is easy to compute but is not required for our purposes).
Applying large deviations bound \cite{large_deviations} to the
i.i.d. sequence $T_j$ we obtain (\ref{eq:LDonT}).
\end{proof}

\bibliographystyle{amsplain}

\begin{thebibliography}{10}





\bibitem{ChenYaoBook}
H.~Chen and D.~Yao, \emph{Fundamentals of queueing networks: Performance,
  asymptotics and optimization}, Springer-Verlag, 2001.

\bibitem{dai}
J.~G. Dai, \emph{On the positive {H}arris recurrence for multiclass queueing
  networks: A unified approach via fluid models}, Ann. Appl. Probab. \textbf{5}
  (1995), 49--77.

\bibitem{DaiMeyn96}
J.~G. Dai and S.~P. Meyn, \emph{Stability and convergence of moments for
  multiclass queueing networks via fluid limit models}, IEEE Transcation on
  Automatic Controls \textbf{40} (1995), 1889--1904.

\bibitem{douformousou04}
R.~Douc, G.~Fort, E.~Moulines, and P.~Soulier, \emph{Practical drift conditions
  for subgeometric rates of convergence}, Adv. Appl. Probab. \textbf{14}
  (2004), no.~3, 1353--1377.

\bibitem{dowmeytwe95a}
D.~Down, S.~P. Meyn, and R.~L. Tweedie.
\newblock Exponential and uniform ergodicity of {M}arkov processes.
\newblock {\em Ann. Probab.}, 23(4):1671--1691, 1995.

\bibitem{dufoco95}
N.~G. Duffield and Neil O'Connell.
\newblock Large deviations and overflow probabilities for the general
  single-server queue, with applications.
\newblock {\em Math. Proc. Cambridge Philos. Soc.}, 118(2):363--374, 1995.

\bibitem{GamarnikZeevi}
D.~Gamarnik and A.~Zeevi, \emph{Validity of heavy traffic steady-state
  approximations in open queueing networks}, {\rm To appear in} Ann. Appl.
  Prob.

\bibitem{ganocowis04}
A.~Ganesh, N.~O'Connell, and D.~Wischik.
\newblock {\em Big queues}, volume 1838 of {\em Lecture Notes in Mathematics}.
\newblock Springer-Verlag, Berlin, 2004.

\bibitem{glywhi94b}
Peter~W. Glynn and Ward Whitt.
\newblock Logarithmic asymptotics for steady-state tail probabilities in a
  single-server queue.
\newblock {\em J. Appl. Probab.}, 31A:131--156, 1994.
\newblock Studies in applied probability.


\bibitem{kumsei90a}
P.~R. Kumar and T.~I. Seidman.
\newblock Dynamic instabilities and stabilization methods in distributed
  real-time scheduling of manufacturing systems.
\newblock {\em IEEE Trans. Automat. Control}, AC-35(3):289--298, March 1990.

\bibitem{kummey96a}
P.~R. Kumar and S.~P. Meyn, \emph{Duality and linear programs for stability and
  performance analysis queueing networks and scheduling policies}, IEEE Trans.
  Automat. Control \textbf{41} (1996), no.~1, 4--17.

\bibitem{malmen79}
V.~A. Maly{\v{s}}ev and M.~V. Men{\cprime}{\v{s}}ikov, \emph{Ergodicity,
  continuity and analyticity of countable {M}arkov chains}, Trudy Moskov. Mat.
  Obshch. \textbf{39} (1979), 3--48, 235, \textit{Trans. Moscow Math. Soc.},
  pp. 1-48, 1981.

\bibitem{CTCN}
S.~P. Meyn.
\newblock Control techniques for complex networks.
\newblock To appear, {Cambridge University Press}, 2007.

\bibitem{mey01a}
S.~P. Meyn, \emph{Sequencing and routing in multiclass queueing networks. {Part
  I}: Feedback regulation}, SIAM J. Control Optim. \textbf{40} (2001), no.~3,
  741--776.

\bibitem{mey05b}
S.~P. Meyn, \emph{Workload models for stochastic networks: Value functions and
  performance evaluation}, IEEE Trans. Automat. Control \textbf{50} (2005),
  no.~8, 1106-- 1122.


\bibitem{metwee}
S.~P. Meyn and R.~L. Tweedie, \emph{Markov chains and stochastic stability},
  Springer-Verlag, 1993.

\bibitem{ResnickStochasticProcesses}
S.~Resnick, \emph{Adventures in stochastic processes}, Birkhäuser Boston, Inc.,
  1992.

\bibitem{rs}
A.~Rybko and A.~Stolyar, \emph{On the ergodicity of stochastic processes
  describing open queueing networks}, Problemi Peredachi Informatsii
  \textbf{28} (1992), 3--26.

\bibitem{large_deviations}
A.~Shwartz and A.~Weiss, \emph{Large deviations for performance analysis},
  Chapman and Hall, 1995.

\bibitem{stolyar}
A.~Stolyar, \emph{On the stability of multiclass queueing networks: A relaxed
  sufficient condition via limiting fluid processe}, Markov Processes and
  Related Fields (1995), 491--512.


\bibitem{tuotwe94a}
P.~Tuominen and R.~L. Tweedie, \emph{Subgeometric rates of convergence of
  $f$-ergodic {Markov} chains}, Adv. Appl. Probab. \textbf{26} (1994),
  775--798.

\end{thebibliography}

\providecommand{\bysame}{\leavevmode\hbox to3em{\hrulefill}\thinspace}
\providecommand{\MR}{\relax\ifhmode\unskip\space\fi MR }
\providecommand{\MRhref}[2]{%
  \href{http://www.ams.org/mathscinet-getitem?mr=#1}{#2}
}
\providecommand{\href}[2]{#2}
\def\cprime{$'$}\def\cprime{$'$}

\end{document}